\documentclass[12pt,letterpaper,reqno]{amsart}

\usepackage{amssymb}
\usepackage{amsmath}
\usepackage{amsthm}
\usepackage{color}
\usepackage{eucal}
\usepackage{graphicx}
\graphicspath{ {pics/} }

\usepackage{caption}
\usepackage{subcaption}

\usepackage{wrapfig}
\usepackage{float}
\usepackage{floatflt}

\usepackage[normalem]{ulem}

\usepackage{mathtools} 

\usepackage{pgfplots}
\pgfplotsset{compat=1.18} 

\newtheorem{theorem}{Theorem}[section]

\theoremstyle{definition}

\newtheorem{remark}[theorem]{Remark}

\title[Cylindrical ZK]{Multi-domain spectral approach\\ for 
Zakharov-Kuznetsov equations in 3D \\ with cylindrical symmetry}

\author[C. Klein]{Christian Klein}
\address{Université Bourgogne Europe, CNRS, IMB UMR 5584, 21000 
Dijon, France,
Institut Universitaire de France} 
\email{Christian.Klein@u-bourgogne.fr}

\author[S. Roudenko]{Svetlana Roudenko}
\address{Department of Mathematics \& Statistics\\
Florida International University,  Miami, FL 33199, USA}
\curraddr{}
\email{sroudenko@fiu.edu}

\author[N. Stoilov]{Nikola Stoilov}
\address{Université Bourgogne Europe, CNRS, IMB UMR 5584, 21000 Dijon, France} 
\email{Nikola.Stoilov@u-bourgogne.fr}



\begin{document}

\begin{abstract}
We present a novel numerical framework for studying nonlinear dispersive equations in higher-dimensional settings, specifically designed for solutions featuring traveling waves along a preferred axis (or field-aligned traveling waves). Using the three-dimensional generalized Zakharov - Kuznetsov (gZK) equation as a model, we convert it into cylindrical coordinates and implement a domain decomposition strategy. 

By partitioning the computational domain into distinct regions based on expected solution behavior, we significantly reduce computational complexity while maintaining the high resolution necessary for capturing small-scale dynamics. Another key innovation of our method is the ability to efficiently handle fractional nonlinearities, specifically, the critical power $p = 7/3$ in 3D, which typically introduces significant computational overhead and numerical instabilities that compromise simulation accuracy. 

Using this framework, we are able to investigate the dynamics of solutions (with cylindrical symmetry) close to the ground state soliton and show that for the 3D critical ZK equation, the ground state serves as the sharp threshold for global vs. finite time existence of solutions. 
Our method successfully tracks the profiles of these singular solutions, providing new insights into the dynamics of wave collapse in three-dimensional magnetized media.

\end{abstract}

\maketitle


\section{Introduction}
This paper proposes a new approach to obtain numerical solutions to the {\it three-dimensional}
Zakharov-Kuznetsov (ZK) equation, 
\begin{equation}\label{ZK3D}
	u_{t}+(\Delta u + u^{p})_{x}=0, \quad (x,y,z) \in \mathbb R^3, ~t \in \mathbb R,
\end{equation}
where $p>1$ is either an integer or a fraction with an odd denominator, $u = u(t,x,y,z)$ is real-valued and $\Delta$ is the standard Laplacian in 3D. 
 
This equation is a three-dimensional variant of the (one-dimensional) generalized  
Korteweg-de Vries (gKdV) equation; 
with the square power, $p=2$, being the celebrated  KdV model for weakly nonlinear waves in shallow water. 
The same square power also appeared in the original 3D ZK equation, which was 
proposed back in 1970's by Zakharov and Kuznetsov in the description of 
weakly magnetized ion-acoustic waves in a low-pressure magnetized 
plasma \cite{ZK1974}, where, in particular, they raised the question of soliton stability  
in a higher-dimensional setting\footnote{The recent work of the second author and collaborators \cite{FHRY} has answered it with the asymptotic stability of solitary waves and the radiation propagating in the opposite direction from the soliton in a cone of $60^o$ angle.}.   
The derivation as a long-wave small-amplitude limit of the Euler-Poisson system in the cold-plasma approximation was done by Lannes, Linares and Saut in \cite{LLS}. 
While the standard ZK equation is characterized by a quadratic nonlinearity ($p=2$ in \eqref{ZK3D}), its physical relevance extends into various other regimes, where the specific particle distribution functions -- such as the presence of superthermal tails, dust grains, or trapped electrons — dictate the power of the nonlinear term. For example, in \cite{V2002} in the unified derivation for multi-species plasmas in 3D, the coefficient of the quadratic term depends strictly on the plasma parameters, such as the density and temperature ratios of the constituent ions and electrons, and thus, may vanish, giving rise to {\it higher} powers in nonlinearity, for example, the modified ZK (mZK) equation with a cubic nonlinearity ($p=3$). 
The inclusion of vortex-like electron distributions leads to the 3D Schamel-type ZK equation, featuring a fractional nonlinearity ($p=\frac32$) that accounts for trapped particle effects \cite{DBD2007b}. Consequently, existence and stability of solitary waves in such models have been questioned in plasma physics literature, e.g., \cite{DBD2007a, DBD2014}. 
Mathematically, these variations are of profound importance in three dimensions, as the power of the nonlinearity determines the stability of solitary waves or instability and potential blow-up, especially as it varies from $p=\frac32$ to $p=3$, since due to scaling as we discuss below, the equation goes from the $L^2$-subcritical to $L^2$-supercritical regime, which affects significantly the global behavior of solutions. 
In this paper we are specifically interested in the 3D $L^2$-critical power, which physically represents the precise threshold where the focusing effects of the plasma medium balance the dispersive decay, making the analysis of these models essential for predicting wave collapse and energy localization in both laboratory, space and astrophysical magnetized plasmas.

Since it is known that the solitary waves of the equation \eqref{ZK3D} travel in the $x$-direction, in this work we consider the 3D ZK equation in {\it cylindrical} coordinates. 
We can introduce in a standard way polar coordinates in 
the $yz$-plane
$$y=\rho\cos\phi,\quad z= \rho\sin\phi.$$ In the case of a 
cylindrical symmetry around the $x$-axis, i.e., $u=u(x,\rho)$, there is no 
dependence on the azimuthal coordinate $\phi$, and hence, the equation \eqref{ZK3D} becomes
\begin{equation}\label{ZKcyl}
u_{t}+\Big(u_{xx}+u_{\rho\rho}+\frac{1}{\rho}u_{\rho}+u^{p}\Big)_{x}=0.
\end{equation}

Equation \eqref{ZKcyl} (or \eqref{ZK3D}) has several conserved quantities, in particular, the $L^2$-norm (or mass),
\begin{equation}
	M[u(t)] = \|u(t)\|_{L^2(\mathbb R^3)}^{2}=2\pi 
	\int_{\rho=0}^{\infty}\int_{x\in\mathbb{R}}^{} \big(u(t, x,\rho))^{2} \, \rho d\rho \, dx \equiv M[u(0)],
	\label{mass}
\end{equation}
and the energy (or Hamiltonian) $E$,
\begin{equation}
	E[u(t)] = 2\pi 
	\int_{\rho=0}^{\infty}\int_{x\in\mathbb{R}}^{} \Big( \frac12(u_{x}^2 + u_{\rho}^2) - \frac1{p+1}u^{p+1} \Big) \, \rho d\rho \, dx
	\equiv E[u(0)].
	\label{energy}
\end{equation}

Among other invariances, equation \eqref{ZK3D} has a scaling invariance (which is also available in cylindrical coordinates): 
if $u(t,x,\rho)$ is a solution of \eqref{ZK3D}, then so is the rescaled version 
\begin{equation}\label{E:scaling}
u_\lambda(t,x,\rho)=\lambda^{\frac32} u(\lambda^3 t, \lambda x, \lambda \rho), \quad \lambda > 0.
\end{equation}
This symmetry makes the Sobolev norm $\dot{H}^s$ with $s=\frac32-\frac2{p-1}$ invariant, 
thus, making the equation \eqref{ZK3D} $\dot{H}^s$-critical. We note that it
is $L^2$-critical ($s=0$) when $p=\frac73$, and that is exactly the nonlinearity we investigate here.

In this paper, we study the $L^{2}$-critical case of the 3D ZK equation numerically, mentioning that the power of nonlinearity is fractional, which for both numerical and analytical purposes makes analysis and computations more challenging.
In our paper \cite{KRS2} we numerically investigated the original 3D {\it quadratic} ZK equation, which is subcritical, and thus, all solutions are global; the power of nonlinearity in that case was even integer, and thus, did not cause any numerical problems. In \cite{KRS2} (see also the 2D work in \cite{KRS1}) we also confirmed numerically the soliton resolution conjecture, showing the splitting of any considered solution into a soliton (or sum of solitons) and outgoing radiation. Analytically, the asymptotic stability of solitons was proved in the energy space $H^1$ for this equation in \cite{FHRY}, and an analytical investigation of the spectral properties for virial operator (with confirming some signs of certain inner products numerically) is done in \cite{HR2026}, see related works in 2D \cite{KRS1}, \cite{CMPS}, \cite{FHRY1}, \cite{HR2025}. 
Thus, by now it is understood how solutions behave in the subcritical case in 3D (including the original 3D quadratic ZK equation), however, as far as the $L^2$-critical ZK equation in 3D, there is not much work, except for an analytical result on certain well-posedness in \cite{LR2021}, and therefore, it is timely to investigate the $L^2$-critical case in 3D, especially, as it is a non-integer power in higher dimensions. This paper is a contribution to such investigations.  
\smallskip

The ZK equation is known to have localized traveling solitary waves of the 
form 
\begin{equation}\label{sol}
u(x,y,z)=Q_c(x-ct,y,z) \equiv Q_c(x-ct, \rho), \quad c>0,
\end{equation}
where $Q_c$ is the rescaling $Q_c(\cdot) = c Q(c \, \cdot)$ with $Q$ a ground state solution of the nonlinear elliptic equation
\begin{equation}\label{E:Q}
\Delta Q - c\,Q+Q^{p}=0.
\end{equation}
By the ground state solution here we mean the unique radial positive $H^1(\mathbb R^3)$ solution of \eqref{sol}, we note that it coincides with the ground state for the 3D nonlinear focusing Schr\"odinger and Klein-Gordon equations. 
Besides radiality and positivity, the properties of $Q$ include $Q \in C^\infty(\mathbb R^3), \partial_r Q(r) < 0$ for any $r = |(x,y,z)|>0$, and for a multi-index $\alpha$ there exists $\gamma>0$ such that
\begin{equation}\label{prop-Q}
|\partial^\alpha Q(x, y, z)| \lesssim_\alpha e^{-\gamma \, r}. 
\end{equation}
From the equation \eqref{E:Q}, it follows that the energy of the ground state is 
\begin{equation}\label{E:E}
E[Q] = \frac{3p-7}{2(5-p)}M[Q],
\end{equation}
and hence, in the $L^2$-critical case ($p=\frac73$), we have $E[Q]=0$, which we confirm numerically in Section \ref{S:solitons}.  
Furthermore, we also obtain numerically that in this critical case 
\begin{equation}\label{E:L2}
\|Q\|_{L^2(\mathbb R^3)} \approx 7.98, 
\end{equation}
consistent with the literature (e.g., see \cite{HR2026}). 

We are specifically interested in the critical case $p=\frac73$, 
where it is expected that the solitary wave solutions \eqref{sol} are 
{\it not stable}. To answer that question numerically, a scheme in {\it higher dimensions}, 
which is accurate, efficient and with high resolution (to track the 
blow-up profile) is needed, and this is exactly what we introduce in this paper. 
\smallskip

The main contribution of this paper is a numerical scheme for a nonlinear 
dispersive equation, which has a {\it cylindrical} symmetry (one preferred 
axis and in addition a dependence only on the distance from this 
axis), which after converting to cylindrical coordinates, splits the 
computational domain into several parts, depending on the expected 
behavior of solutions and then gluing them at least at the $C^1$ 
level. We use different independent variables in these domains for an 
optimal allocation of numerical resolution and to obtain less 
singular forms of the equation on the symmetry axis. 
This allows us to efficiently construct the time evolution of solitary waves and their perturbations, thus, answering the question about the stability of the solitons, for example, in the 3D critical ZK equation. Furthermore, we also check nearby generic (Gaussian) data to further confirm stability/instability regions and long-term behavior of cylindrically symmetric solutions, finding that indeed if 
\begin{equation}\label{E:massQ}
\|u_0\|_{L^2(dx,d\rho)} \leq \|Q\|_{L^2(\mathbb R^3)},
\end{equation}
then solutions exist for all time and disperse (or scatter). If the 
inequality in \eqref{E:massQ} is of the opposite sign, then solutions that we considered with the initial data $u_0$ 
(namely, perturbations of ground state and Gaussian) numerically blow up in finite time in $L^{\infty}$. 
\smallskip

The paper is organized as follows: in Section \ref{S:Num} we present our new numerical approach used to solve the 3D ZK equation with cylindrical symmetry. 
In Section \ref{S:solitons} we first show the construction of the 
ground state and its evolution, and then in Section \ref{S:sol-pert} 
consider several perturbations. Following that, in Section 
\ref{S:gaussian} we include examples of gaussian initial data, to 
confirm the ground state threshold for this type of data.   We add 
some concluding remarks in section 6. 
\smallskip

{\bf Acknowledgments.} The work of C. K. and N.S. was partially supported by the ANR project 
ANR-17-EURE-0002 EIPHI and by the ANR project 
ISAAC-ANR-23-CE40-0015-01. 
S.R. would like to thank the hospitality and support of IMB, where most of this work was done. S.R. was partially supported by the NSF grant DMS-2452782.

\section{Numerical approach}\label{S:Num}
In this section, we introduce
the numerical approach using a 
Fourier spectral method in the $x$ variable and a {\it multi-domain} splitting 
for the 
$\rho$-dependence together with an implicit fourth order Runge-Kutta 
method for the time integration. 

\subsection{Fourier spectral approach}
The dependence of the solution on the coordinate $x\in\mathbb{R}$ is 
approximated in standard way on a torus $\mathbb{T}_{2\pi L}$, where 
$L>0$ is chosen large enough that the studied solution can be treated 
as essentially periodic for localized functions as the solitary wave. 
On the torus we introduce the standard discretization of the discrete 
Fourier transform that will be efficiently computed with a Fast 
Fourier Transform (FFT): $x_{n}=L(-\pi+2\pi n/N)$, $n=1,\ldots,N$, 
$N\in\mathbb{N}$. As it is well-known the FFT coefficients denoted by 
$\hat{u}_{n}$ will decrease exponentially with $n$ for smooth periodic 
functions $u(x)$. This allows us to estimate the numerical error via the 
highest FFT coefficients, see for instance the discussion in 
\cite{trefethen}. 

\subsection{Chebyshev collocation method}

The appearance of the term $\frac1{\rho}u_{\rho}$ in \eqref{ZKcyl}, i.e., of 
a singular and non-autonomous term in the equation, makes the use of 
Fourier spectral methods in this variable inefficient. Instead, we 
apply a Chebyshev collocation method to treat the dependence in 
$\rho$, see \cite{trefethen} and references therein for details. The 
idea is to introduce for a variable $l\in[-1,1]$ the collocation 
points 
$$l_{n}=\cos(\pi n/N_{c}),\quad n=0,1,\ldots,N_{c},$$
and to consider the Lagrange interpolation polynomial for a function 
$f(l)$ on these collocation points. The derivative of the Lagrange 
polynomial is used to approximate the derivative of $f$ with respect 
to $l$ which leads to the standard \emph{Chebyshev differentiation 
matrices}, see \cite{trefethen,WR}. This approach is equivalent to 
the expansion of the function $f$ in terms of Chebyshev polynomials 
$T_{n}(l)$,
$f(l)\approx \sum_{n=0}^{N_{c}}a_{n}T_{n}(l)$, where 
$T_{n}(l)=\cos(n\arccos(l))$. The Chebyshev coefficients $a_{n}$ can 
be determined with a collocation method, i.e., by imposing on the 
collocation points
$$f(l_{m})=\sum_{n=0}^{N_{c}}a_{n}T_{n}(l_{m}),\quad 
m=0,\ldots,N_{c}.$$ 
The coefficients can be conveniently obtained with a 
\emph{Fast cosine transform} (FCT) related to the FFT. It is known 
that the Chebyshev coefficients $a_{n}$ of an analytic function 
decrease exponentially with $n$, see the discussion in 
\cite{trefethen}. Thus, the coefficients of order $N_{c}$ indicate as 
the FFT coefficients the order of the numerical error.

\subsection{Multi-domain spectral approach}
To efficiently approximate the $\rho$-dependence and to take care of 
the singular term, we introduce two domains in $\rho$ and  
different independent variables in each domain:

\textbf{Domain I:} $\rho<\rho_{0}$ (we generally work with 
$\rho_{0}=1$ here).

In this domain we use, as in \cite{CKS}, the variable $s=\rho^{2}$, in 
which equation (\ref{ZKcyl}) reads 
\begin{equation}\label{ZKs}
	u_{t}+(u_{xx}+4su_{ss}+4u_{s}+u^{p})_x=0,
\end{equation}
a less singular form of the equation. Since it is known that 
functions regular in the vicinity of the symmetry axis are analytic in 
$\rho^{2}$, this also leads to an optimised allocation of numerical 
resolution. We use the 
transform $s=\rho_{0}^{2}(1+l)/2$, with $l\in[-1,1]$, and introduce 
$N_{I}+1$ Chebyshev collocation points for $l$. The solution in this domain is then denoted by $u^I$.

\textbf{Domain II:} $\rho_{0}<\rho<\rho_{1}$ ($\rho_{1}\gg1$). 

In this domain $\rho$ does not vanish, and the term with $1/\rho$ 
does not pose numerical problems if $\rho_{0}\sim 1$. 
We apply the relation $\rho=\rho_{0}(1+l)/2+\rho_{1}(1-l)/2$ with 
$l\in[-1,1]$ and introduce $N_{II}+1$ Chebyshev collocation points 
for $l$. The solution in this domain is denoted by $u^{II}$.

Since the equation is singular only for $s=0$, no condition needs to be 
given there. However, a condition needs to be imposed for 
$\rho=\rho_{1}$. We generally take $\rho_{1}$ sufficiently large that 
the solution can be assumed to vanish there with numerical accuracy. 
In cases where this is not possible, one could use a compactified 
third domain as in \cite{BK,KS2021}, but this was not necessary for the 
problems studied in this paper. 

For $\rho=\rho_{0}$ we impose that the solution is $C^{1}$ in $\rho$. 
This means we impose the conditions
\begin{equation}
	u^{I}(\rho_{0}^{2})=u^{II}(\rho_{0}) \quad \mbox{and} \quad 
	2\rho_{0}u^{I}_{s}(\rho_{0}^{2})=u^{II}_{\rho}(\rho_{0}).
	\label{C1}
\end{equation}
These two conditions as well as the vanishing condition for 
$\rho_{1}$ are imposed via a $\tau$-method (see for instance 
\cite{trefethen}): lines in the differentiation matrices approximating 
$\partial_{\rho\rho}+\frac1{\rho}\partial_{\rho}$ and 
$4s\partial_{ss}+4\partial_{s}$ are replaced by these conditions. 
This assures that the boundary and matching conditions are 
approximated with the same accuracy as the differential equation. 

Note that Chebyshev differentiation matrices have a conditioning of 
the order of $\mathcal{O}((N_{c})^{2})$, which implies that the value 
of $N_{c}$ should be kept small in order to avoid a loss of accuracy 
due to this conditioning. This can be addressed by introducing more 
domains as above, but was not needed here.

\subsection{Time integration}
The spatial discretization introduced above leads to a finite 
dimensional system of ordinary differential equations in time. Due to 
the third derivative in $x$ and the use of Chebyshev differentiation 
matrices for $\rho$, this leads to a \emph{stiff} system, which loosely 
speaking means that the explicit time integration schemes will be 
inefficient due to stability conditions. Since we are also interested in 
studying blow-up phenomena, we apply a 4th order method, an 
implicit Runge--Kutta scheme of order 4 
(IRK4).  The method is also known as Hammer-Hollingsworth or 2-stage 
Gauss method. We introduce $N_{t}$ time steps $t_{n}$, 
$n=0,\ldots,N_{t}$ and denote $h:=t_{n+1}-t_{n}$. 

The general formulation of an $s$-stage Runge--Kutta method for the initial value problem
$y'=f(y,t),\,\,\,\,y(t_0)=y_0$ is the following:
\begin{eqnarray}
 y_{n+1} = y_{n} + h      \underset{i=1}{\overset{s}{\sum}} \, 
 b_{i}K_{i}, \\
 K_{i} = f\left(t_{n}+c_ {i}h,\,y_{n}+h  
 \underset{j=1}{\overset{s}{\sum}} \, a_{ij}K_{j}\right),
 \label{Ki}
\end{eqnarray}
where $b_i,\,a_{ij},\,\,i,j=1,...,s$, are real numbers and
$c_i=   \underset{j=1}{\overset{s}{\sum}} \, a_{ij}$.  
Here, we take $s=2$ and
\begin{equation}
c_{1}=\tfrac{1}{2}-\tfrac{\sqrt{3}}{6}, 
c_{2}=\tfrac{1}{2}+\tfrac{\sqrt{3}}{6}, a_{11}=a_{22}=\tfrac14,
a_{12}=\tfrac{1}{4}-\tfrac{\sqrt{3}}{6}, 
a_{21}=\tfrac{1}{4}+\tfrac{\sqrt{3}}{6}, ~~\mbox{and} ~~ b_{1}=b_{2}=\tfrac12.
\end{equation}
This scheme is of classical order 4, but stage order 2. 

We consider equation \eqref{ZKcyl} after an FFT in $x$ in the form
\begin{equation}\label{ZKhat}
\hat{u}_{t}+ik(-k^{2}+\mathcal{L})\hat{u}=-ik\widehat{\,u^{p}\,},
\end{equation}
where $\mathcal{L}$ stands for the differentiation matrices introduced 
in the domains I and II. Note that equation \eqref{ZKhat} decouples 
in $k$, thus, allowing us to treat each of the values of $k$ separately in 
a parallel computation. 

Equation \eqref{Ki}, for instance, for $K_{1}$ is written in the form 
\begin{equation}\label{step}
(1 - iha_{11}k(\mathcal{L}-k^2))\hat{K}_{1}=-ik\widehat{\,\tilde{u}_{1}^{p}\,} - 
ik(n)(\mathcal{L}-k^2)(\hat{u}_{n}+ha_{12}\hat{K}_{2}),
\end{equation}
where the index denotes the time step and the tilde notation means
$\tilde{u}_{i}=u_{n}+ha_{i1}K_{1}+ha_{i2}K_{2}$, $i=1,2$. An analogous 
equation can be written for $\hat{K}_{2}$. This leads to a nonlinear system 
for $K_{1},K_{2}$, which is then solved iteratively with a simplified 
Newton iteration: the operator on the left-hand side of \eqref{step} 
is inverted in each step of the iteration. The $\tau$-method 
is implemented in exactly this inversion. Thus, the computational 
cost per iteration will consist of inverting two 
$(N_{I}+N_{II}+2)\times(N_{I}+N_{II}+2)$ matrices that consist
except for the matching conditions essentially of blocks of size 
$(N_{I}+1)\times(N_{I}+1)$ and $(N_{II}+1)\times(N_{II}+1)$. The 
iteration is stopped once the change in the $K_{1}, K_{2}$, per iteration 
is smaller than some threshold, typically $10^{-6}$. Note that this 
relatively large residual for the iteration allows to reach much 
higher accuracies than this since there is a factor $h$ in front of 
the $K_{i}$, $i=1,2$ in (\ref{Ki}). The convergence of the simplified 
Newton iteration is in general considerably faster than a fix point 
iteration (if the latter converges at all). 

The numerical error is estimated via the conservation of the 
numerically computed energy, which will depend on time due to 
unavoidable numerical errors. As discussed in \cite{etna,KR}, the 
relative conservation of the energy will in general overestimate the 
actual accuracy by 1-2 orders of magnitude.  

\begin{remark}\label{R:2}
	The FFT and its inverse introduce some {\it small imaginary} part in 
	the transformed real solution. This can lead to problems if  
	fractional nonlinearity is considered, which is our case, since $p=7/3$. To avoid 
	this issue (which can significantly influence the conserved quantities),
	we always enforce real values after an inverse FFT. In 
	order to take care of the Matlab root that is branched at the 
	negative real axis, we compute $u^{7/3}$ as $|u|^{7/3}\mbox{sign}(u)$ 
	and the term $u^{10/3}$ as $|u|^{10/3}$. 
\end{remark}

\section{Solitary wave profiles}\label{S:solitons}

\subsection{Profiles construction}
\begin{figure}[htb!]
\includegraphics[width=0.7\textwidth, height=.48\textwidth]{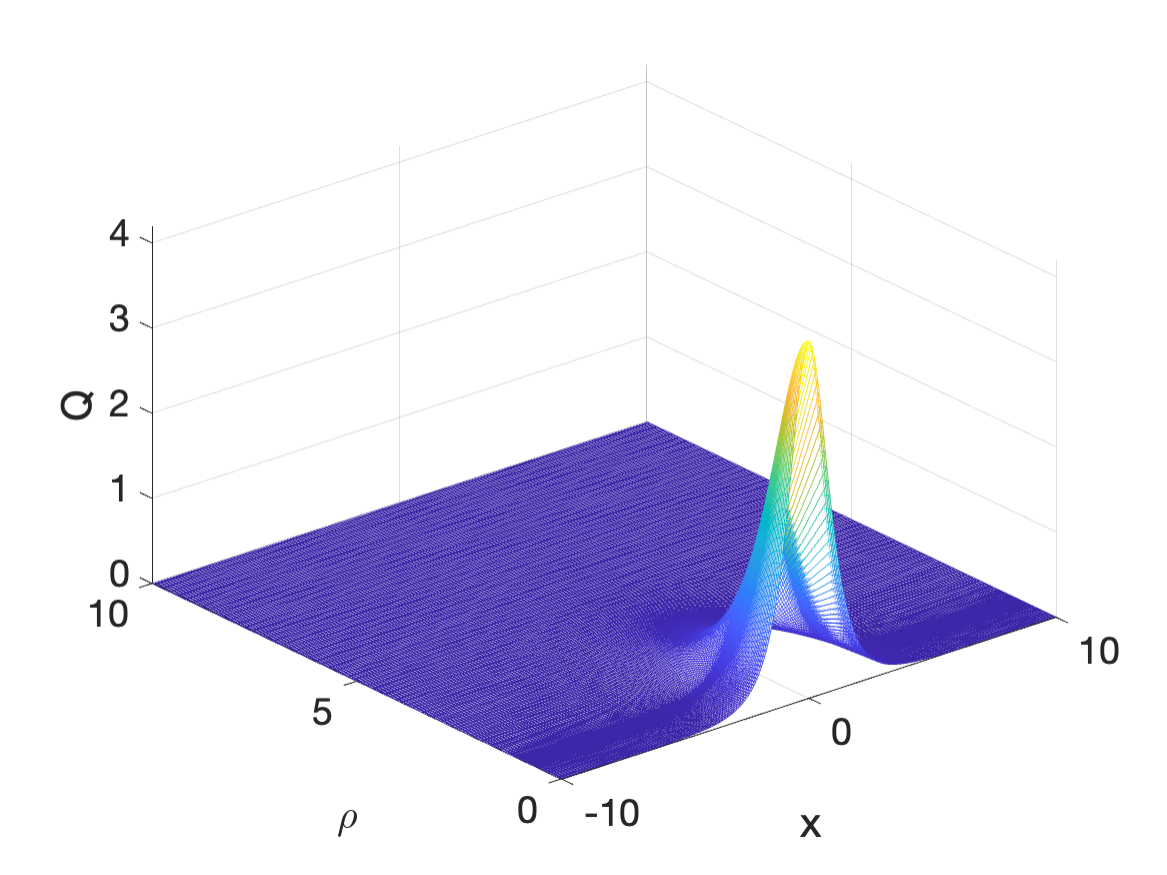}
\caption{Profile of the ground state for the 3D ZK equation \eqref{ZK3D}, $p=\frac73$.}
\label{solfig}
\end{figure}

Since the profiles of solitary waves \eqref{sol} are radially symmetric, they can 
be numerically constructed via different approaches, for instance, via the 
one used in \cite{CKS} or in \cite{RRY} and \cite{RWY}. The radially symmetric solution can be then 
interpolated on the grid of the above numerical approach. 
Alternatively, we can also solve equation \eqref{sol} directly with the 
numerical approach outlined above: the discretization in $x$ and 
$\rho$ leads to a system of nonlinear equations for $Q$, which is 
solved with a Newton-Krylov iteration, where the action of the inverse 
of the Jacobian is iteratively computed with GMRES \cite{GMRES}. We 
use $N=2^{9}$ FFT modes for $x\in 5[-\pi,\pi]$, $N_{I}=20$ 
Chebyshev polynomials for $s<1$ and $N_{II}=100$ Chebyshev 
polynomials for $1<r<20$. We show the resulting solution in 
Fig.~\ref{solfig}. Its mass is $M\approx63.7831$ (which is the square of the value in \eqref{E:L2}; its energy vanishes 
(numerically it is on the order of $10^{-11}$, indicating 
the high accuracy, with which it has been constructed) consistent with \eqref{E:E}.
\begin{figure}[htb!]
\includegraphics[width=0.49\textwidth]{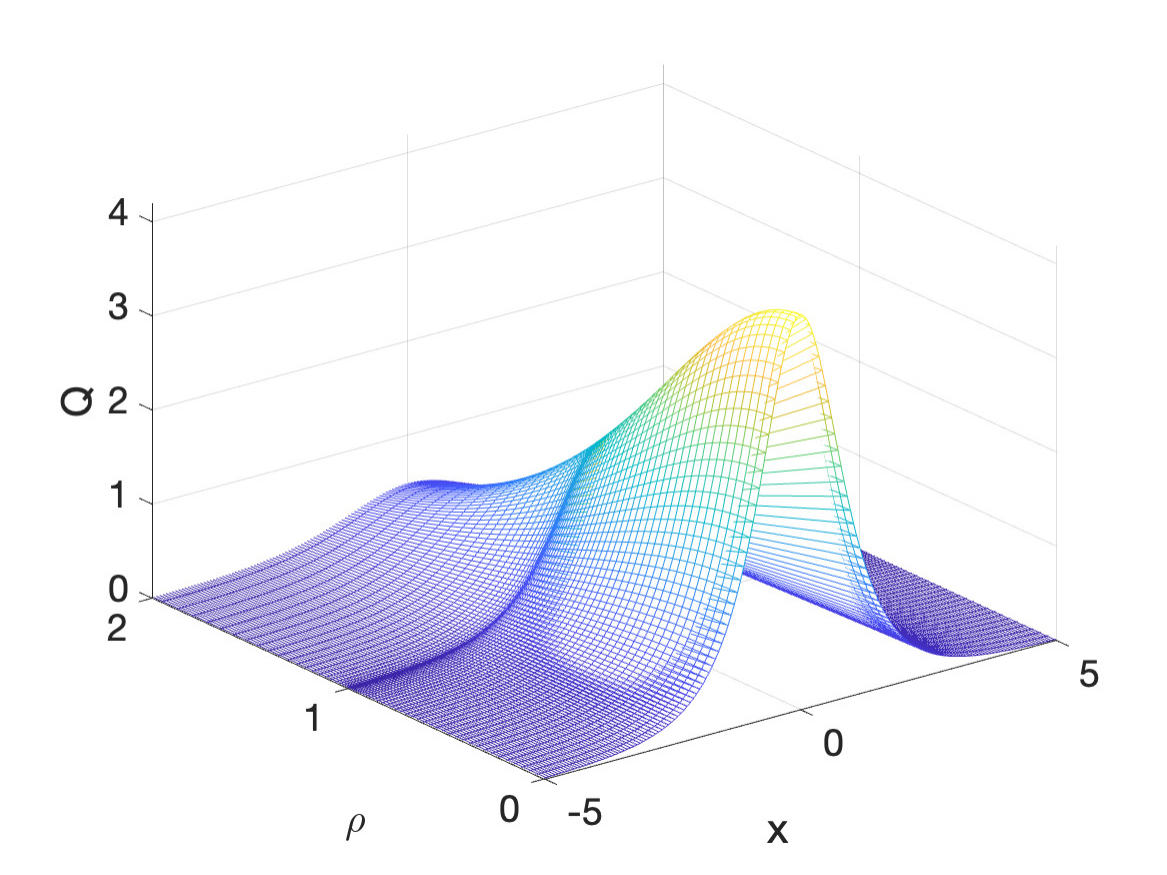}
\includegraphics[width=0.49\textwidth]{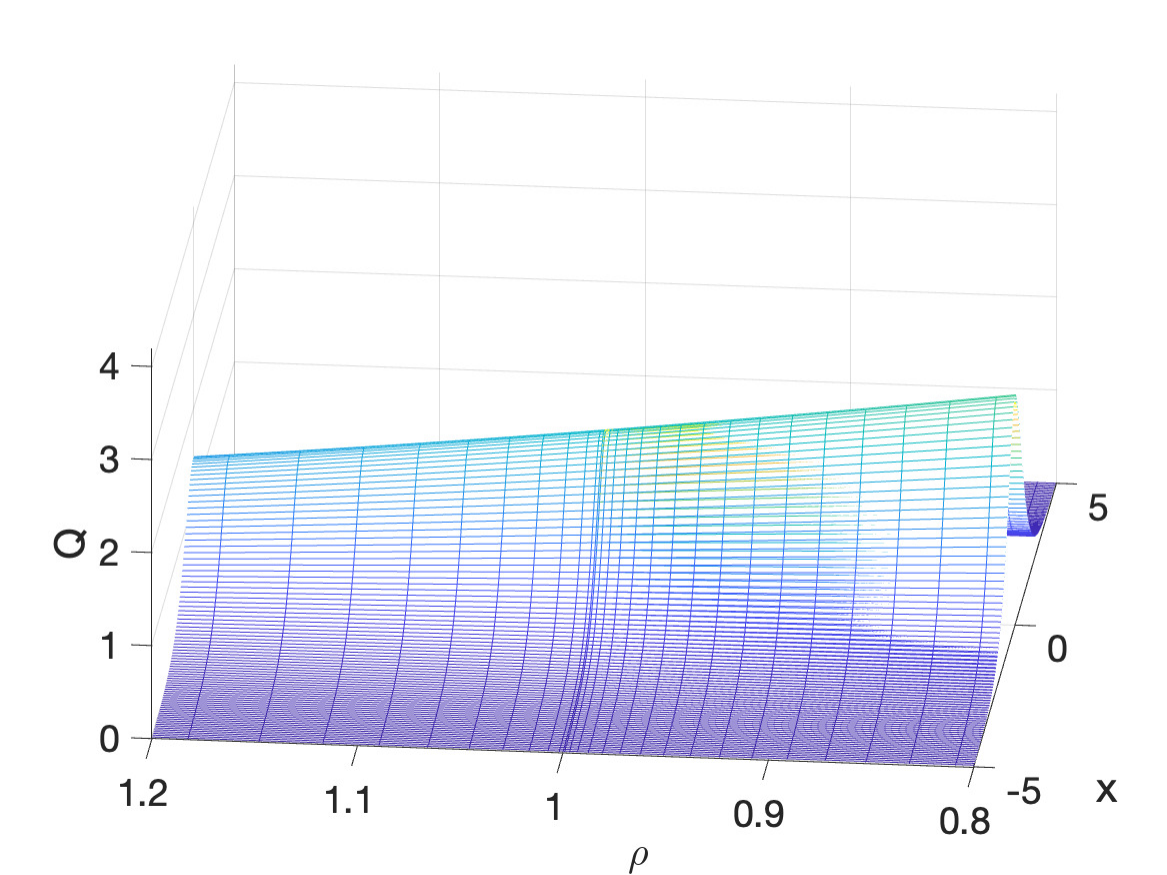}
\caption{Zoom-in into the ground state profile from Fig.~\ref{solfig} for the 3D ZK equation \eqref{ZK3D}, $p=\frac73$.}
\label{solfigcu}
\end{figure}

Close-ups of the solution, illustrating the splitting of the computational domain into subdomains I and II, and the 
boundary between I and II, are shown in Fig.~\ref{solfigcu}. One can 
see that the domain I is set up so that it could be more refined to track the growth of the solution, due to the blow-up behavior. Also, note on the right plot of Fig. \ref{solfigcu} the smoothness of the solution where two domains are attached to each other (at $\rho_0=1$), the solution itself and its derivative are smoothly glued.

\begin{figure}[htb!]
 \includegraphics[width=0.49\textwidth]{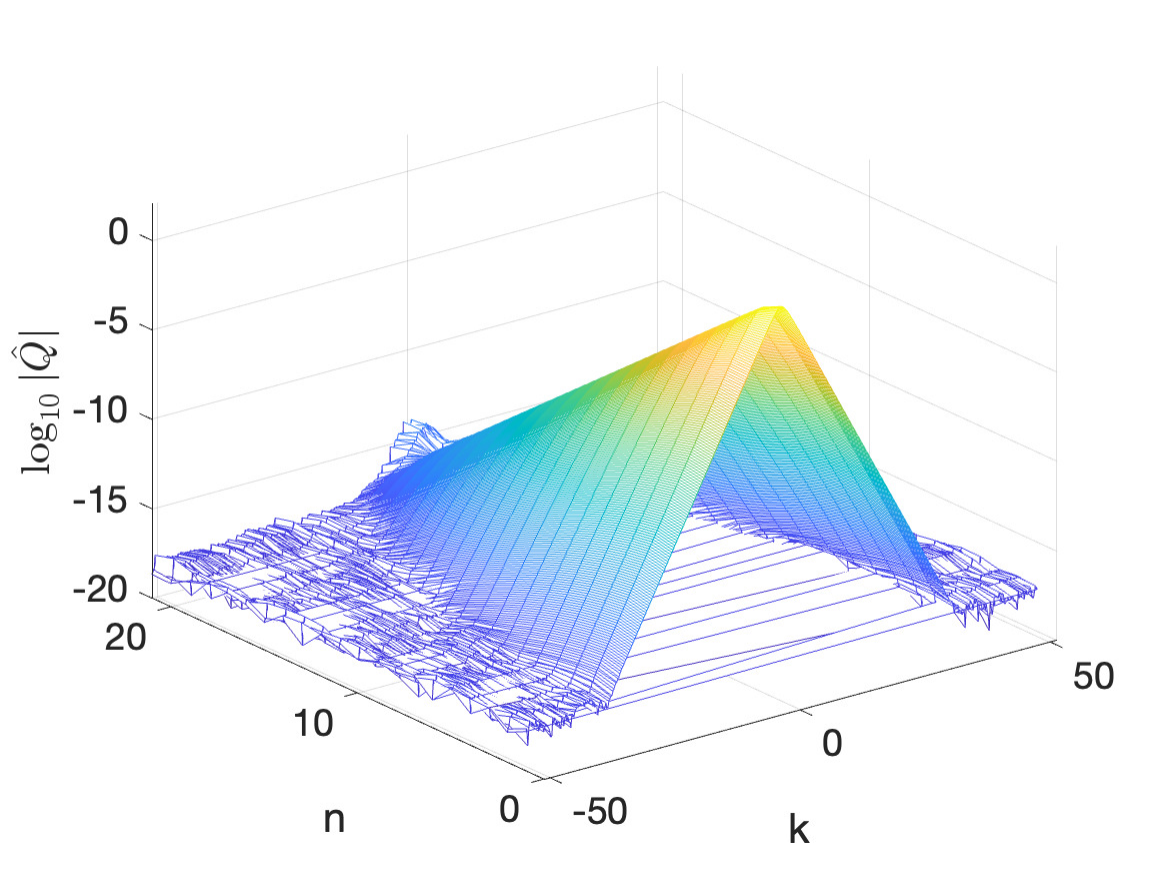}
 \includegraphics[width=0.49\textwidth]{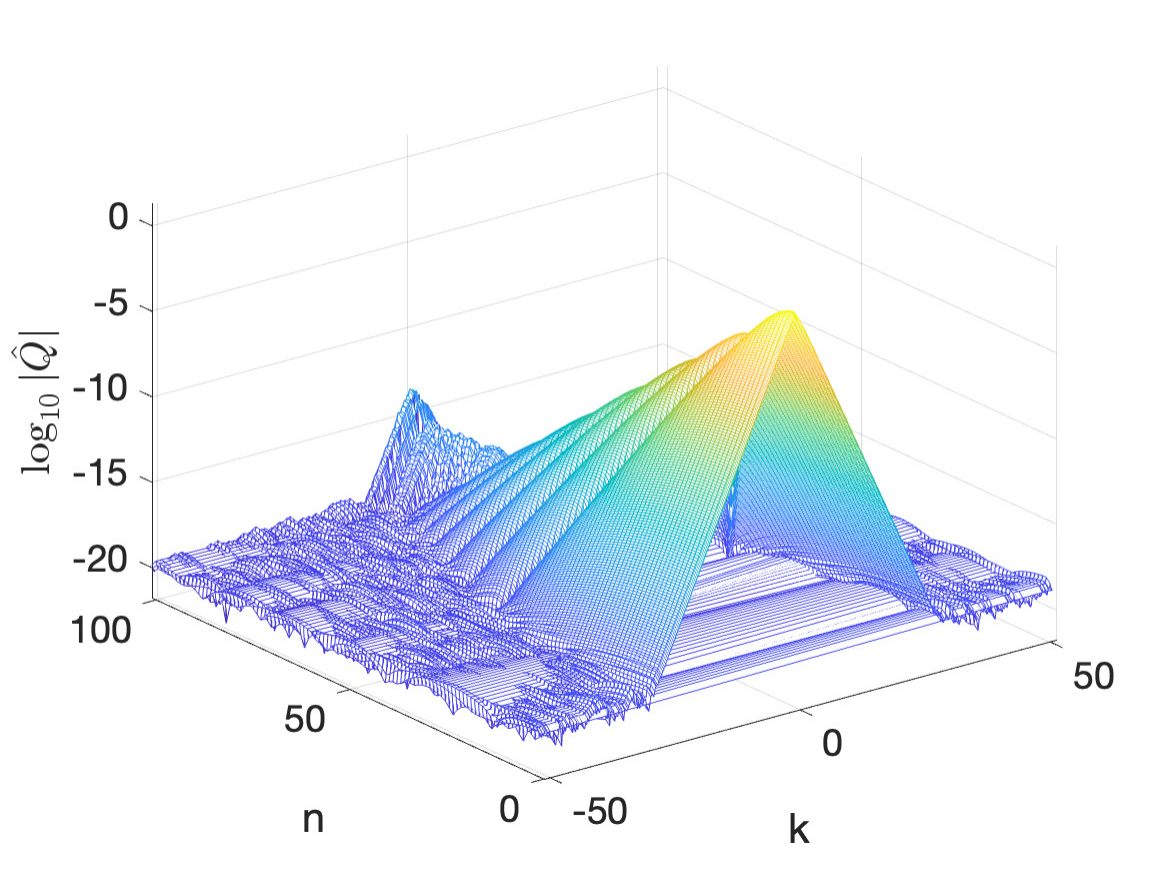}
 \caption{Spectral coefficients of the solitary wave of 
 Fig.~\ref{solfig}, on the left in domain I, on the right in domain II.}
 \label{solfigcoeff}
\end{figure}

The spectral coefficients of the solution are shown in the two 
plots of Fig.~\ref{solfigcoeff}. It can be seen that the 
coefficients decrease to machine precision for both: large $|k|$ and 
$n$, showing that the spatial dependence is fully resolved 
numerically. 

\subsection{Dynamical test of profiles}

\begin{figure}[htb!]
\includegraphics[width=0.65\textwidth]{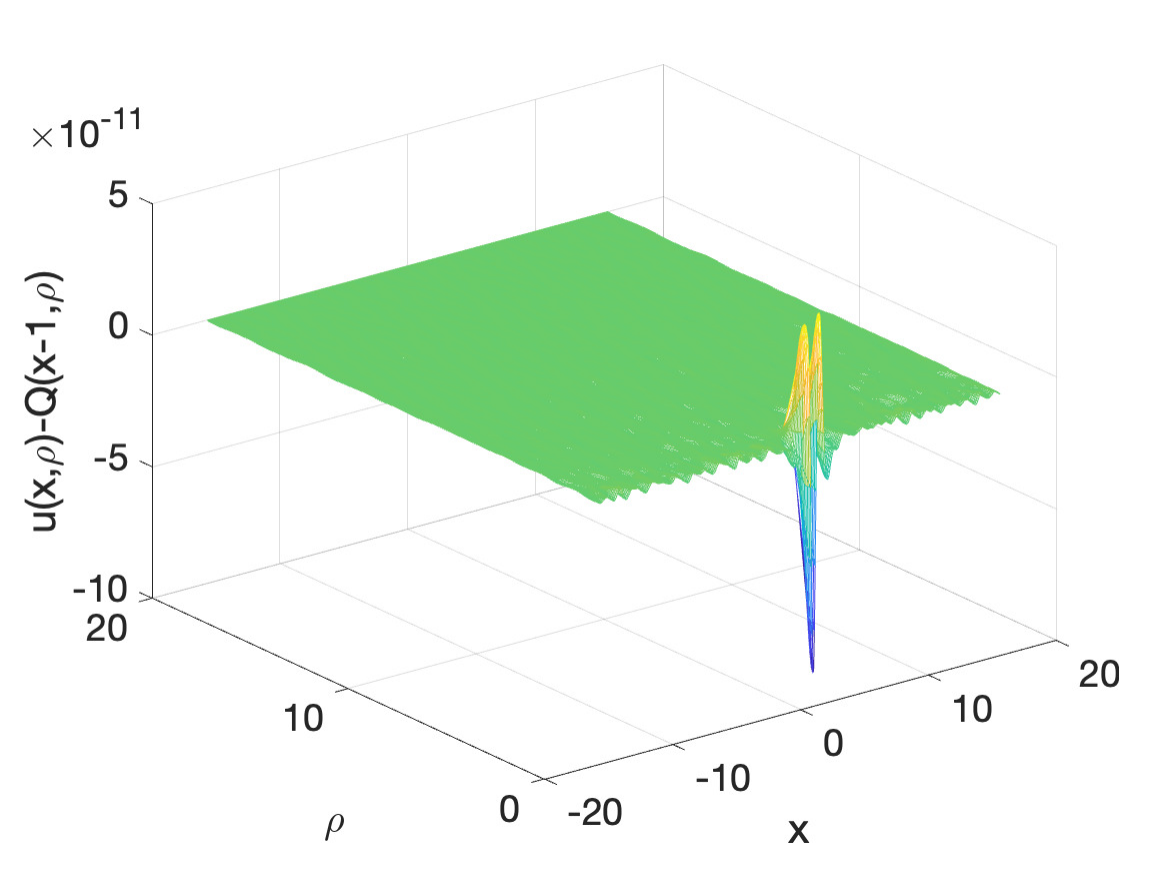}
\caption{Difference of the solution for solitary wave initial data 
at $t=1$ and the shifted solitary wave.}
\label{soldifft1}
\end{figure}

We use the profile in Fig.~\ref{solfig} of the solitary wave (aka ground state) as initial data for the time evolution 
code ($t=0$) and propagate (with the critical ZK) the solution with $N_{t}=400$ time steps to $t=1$. 
Comparing the numerical solution at the final time with the shifted 
(velocity 1) initial data, we find these difference to be on the order of 
$10^{-10}$, see Fig.~\ref{soldifft1}.
This also tests the numerically constructed ground state (or profile of the solitary wave), which 
was assumed to be constructed at least with this accuracy. The energy and the 
mass are conserved to the order of $10^{-12}$. This shows that the 
code is able to study the time evolution of solutions to the 3D ZK 
equation (with a fractional nonlinear power) with a precision of $10^{-10}$ or better.

\section{Perturbations of the solitary wave}\label{S:sol-pert}
In this section we use our method to study amplitude perturbations of the solitary wave, first by 
considering initial data of the form 
$$
u(0, x,\rho)=\lambda Q(x,\rho) \quad \mbox{with ~ real}\quad \lambda\sim1.
$$ 

\begin{figure}[htb!]
 \includegraphics[width=0.5\textwidth]{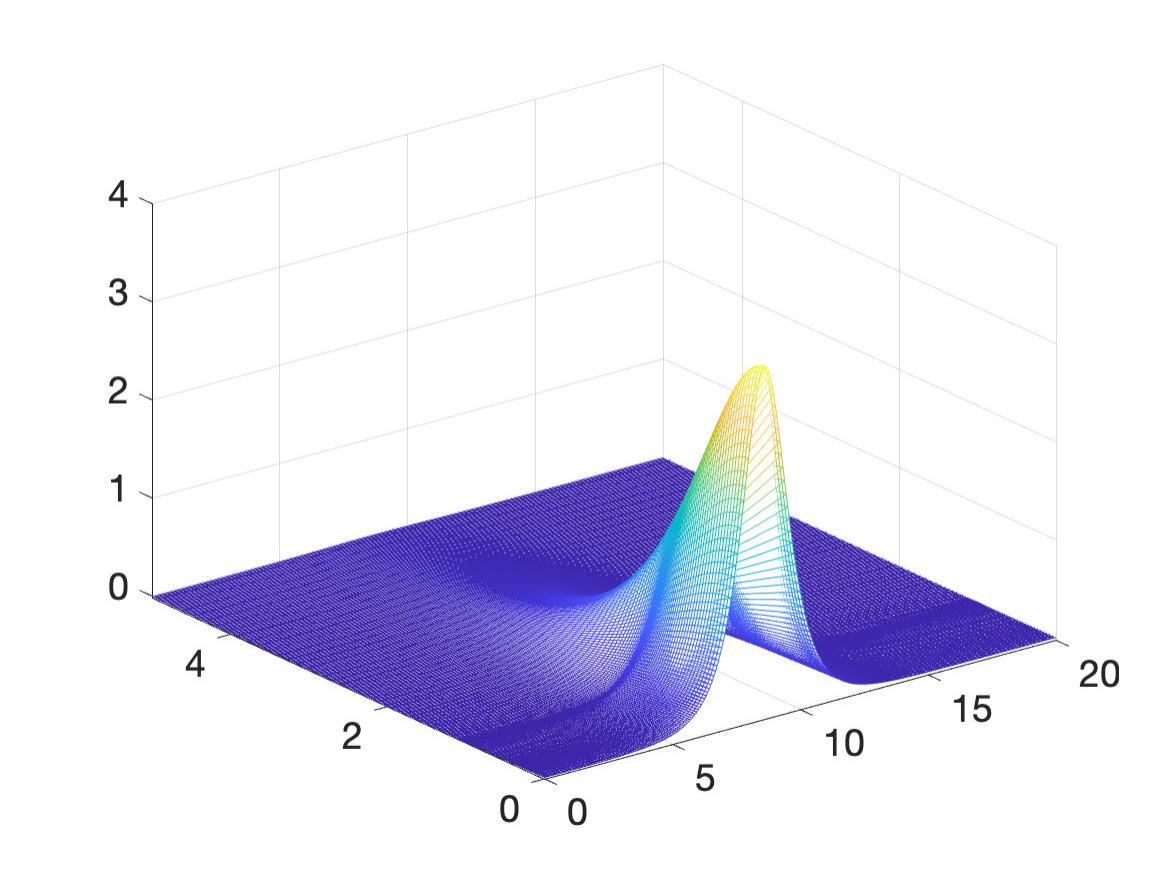}
 \includegraphics[width=0.48\textwidth]{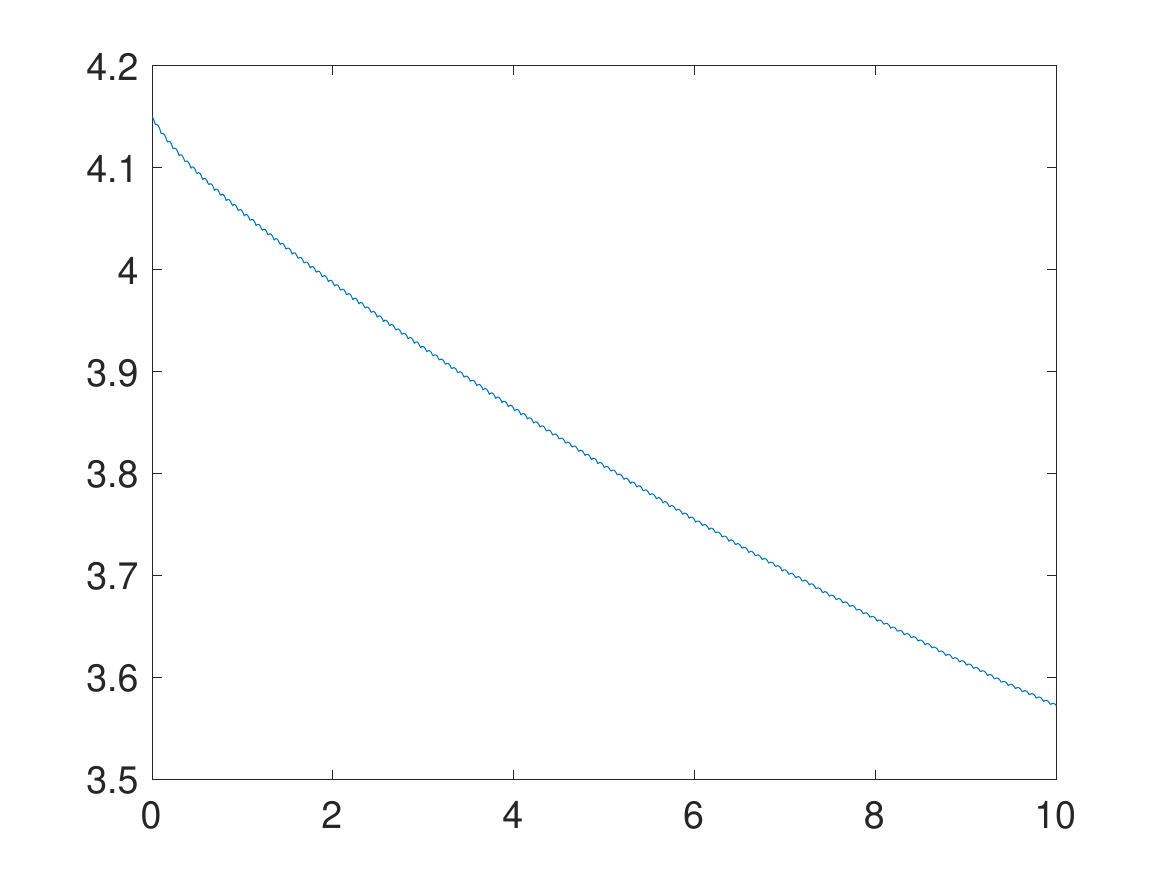}
 \caption{Solution to the 3D critical ZK equation with initial condition 
 $u(0,x,\rho)=0.99\, Q(x,\rho)$, on the left the solution at $t=10$, 
 on the right the time dependence of the $L^{\infty}$ norm.} 
 \label{sol099}
\end{figure}

For the case of $\lambda=0.99$, a perturbation with smaller 
mass than the original ground state mass, 
we use the same numerical 
parameters as before and $N_{t}=1000$ time steps for $t\leq 10$. We 
show a snapshot of the solution of the 3D critical ZK equation for this initial condition 
at $t=10$ in Fig.~\ref{sol099} on the left. The initial condition visibly disperses and 
flattens. This is confirmed by the $L^{\infty}$ norm of the solution 
on the right of the same figure. It appears to decrease 
monotonically, the initially data seem to be simply radiated away to 
infinity. 

\begin{figure}[htb!]
\includegraphics[width=0.5\textwidth]{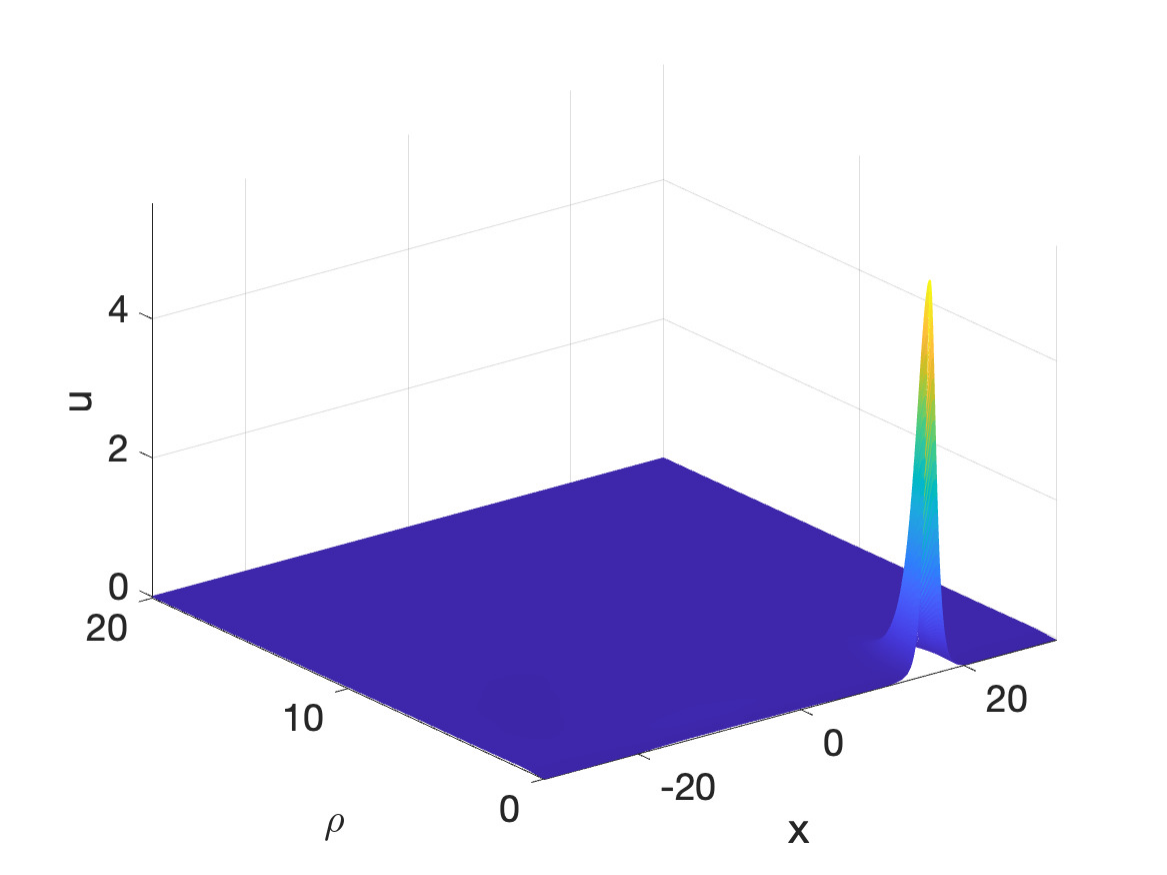}
\includegraphics[width=0.49\textwidth]{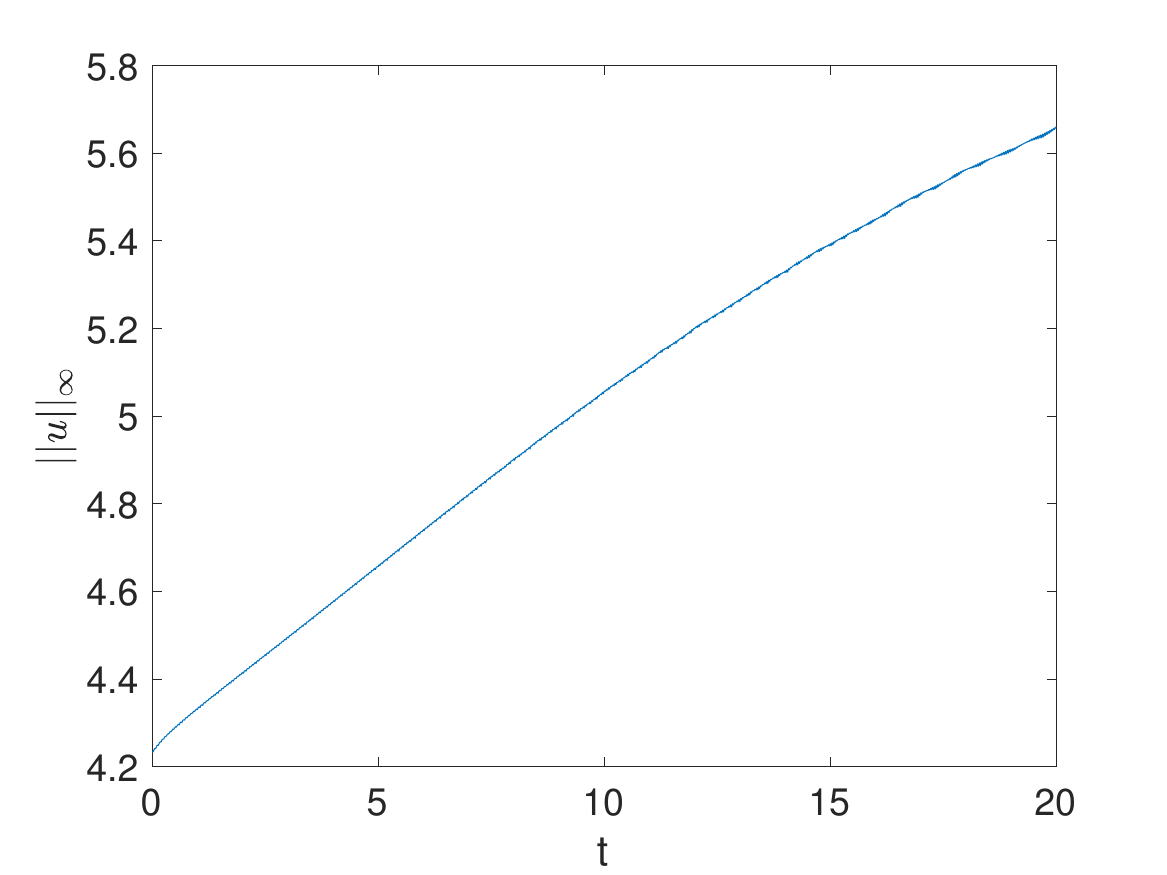}
\caption{Solution to the 3D critical ZK equation with initial condition 
 $u(0,x,\rho)=1.01\, Q(x,\rho)$, on the left the solution at the 
 final computational time $t= 2$, 
 on the right the time dependence of the $L^{\infty}$ norm up to $t=20$.}
\label{sol1.01}
\end{figure}

\begin{figure}[htb!]
\includegraphics[width=0.49\textwidth]{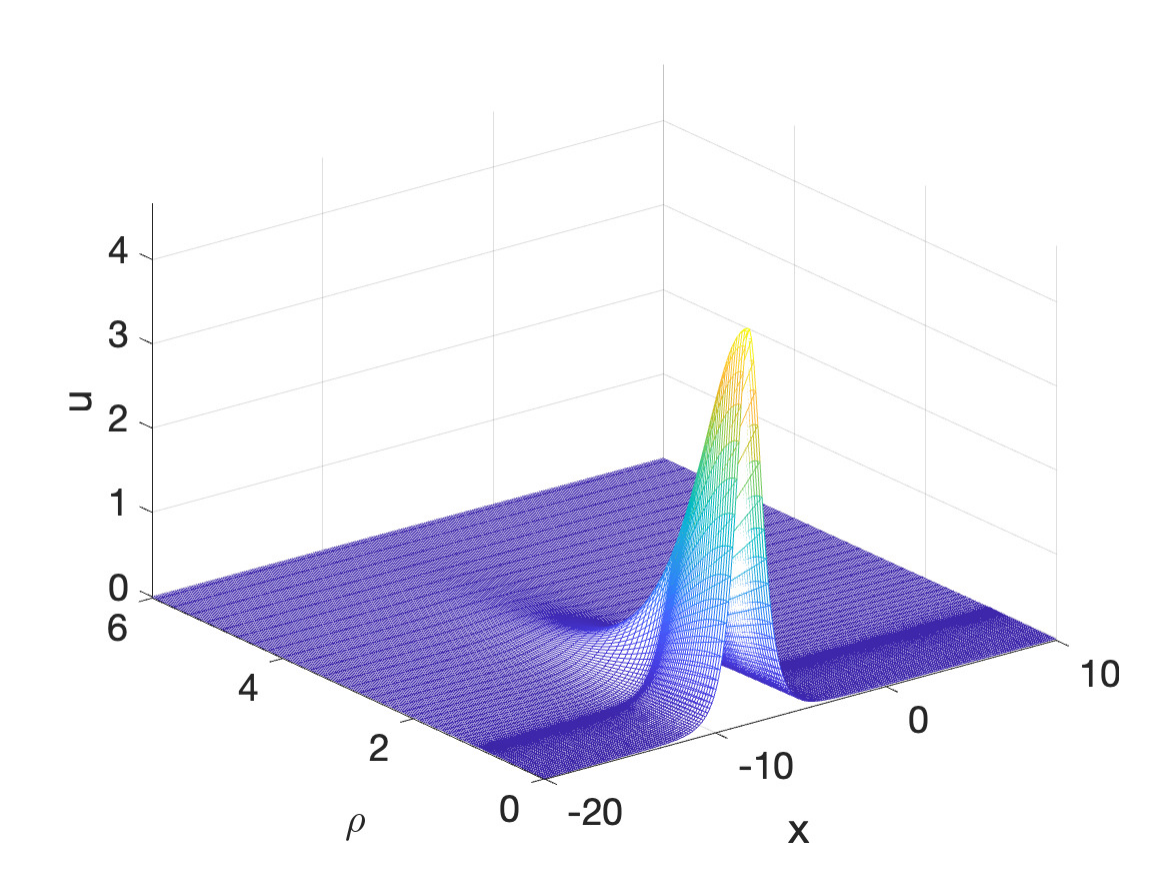}
\includegraphics[width=0.49\textwidth]{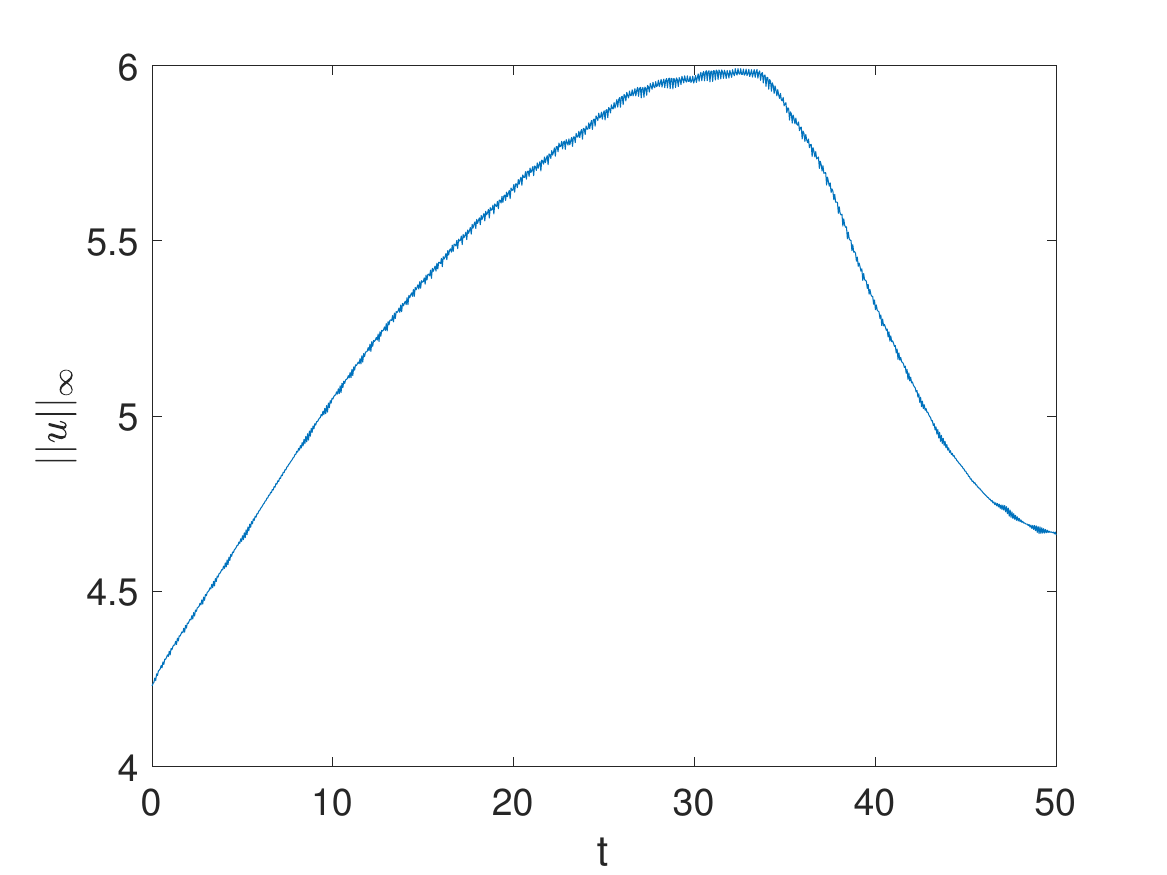}
\caption{Solution to the 3D critical ZK equation with initial condition 
 $u(0,x,\rho)=1.01\, Q(x,\rho)$, on the left the solution at the 
 final computational time $t=50$, 
 on the right the time dependence of the $L^{\infty}$ norm up to $t=50$.}
\label{sol1.01add}
\end{figure}

If we consider a perturbation of the solitary wave with a slightly 
larger mass,  for instance $\lambda=1.01$, the $L^{\infty}$ norm of 
the solution will initially grow, see Fig.~\ref{sol1.01} on the 
right. Since our numerical simulations are done on a finite interval, a similar 
behavior can be observed in this case as for the NLS equation on a finite 
interval, see \cite{KRSfinite} and references therein: perturbations 
of the soliton that blow up on $\mathbb{R}$ are actually stable on a 
finite interval.
A similar phenomenon can be observed here, 
partially due to the radiation re-entering the domain or being 
reflected at $\rho=\rho_{1}$, 
and thus, 
arresting the blow-up for very small mass-supercritical perturbations. 
For longer times, due to the periodicity in $x$, the solution 
leaves the computational domain at $L\pi$ and reenters it at $-L\pi$. 
The solution for $t=20$ is shown on the left of Fig.~\ref{sol1.01add}.  
In the same figure, right plot, it can be seen that the $L^{\infty}$ norm 
reaches a maximum (since the initial mass only slightly larger than 
that of the ground state) and then it starts decreasing. Note that 
the oscillations of the $L^{\infty}$ norm are due to the radiation 
present in the computational domain and due to the fact that the 
maximum is evaluated on the collocation points (the actual maximum of 
the solution might not be on such a point).
Thus, similar to the NLS on bounded domains \cite{KRSfinite}, we conclude that on the 
considered finite computational domain, the 3D critical ZK soliton is actually stable for very
small perturbations, due to the boundedness of the domain (the 
$L^{\infty}$ norm of the solution will oscillate around some mean 
value indicating a stable soliton as in \cite{KRSfinite}).    

To check whether it is actually unstable in the infinite space setting, one has either to consider 
much larger domains or perturbations with larger mass. Below we consider larger mass perturbation.

\begin{figure}[htb!]
 \includegraphics[width=0.5\textwidth]{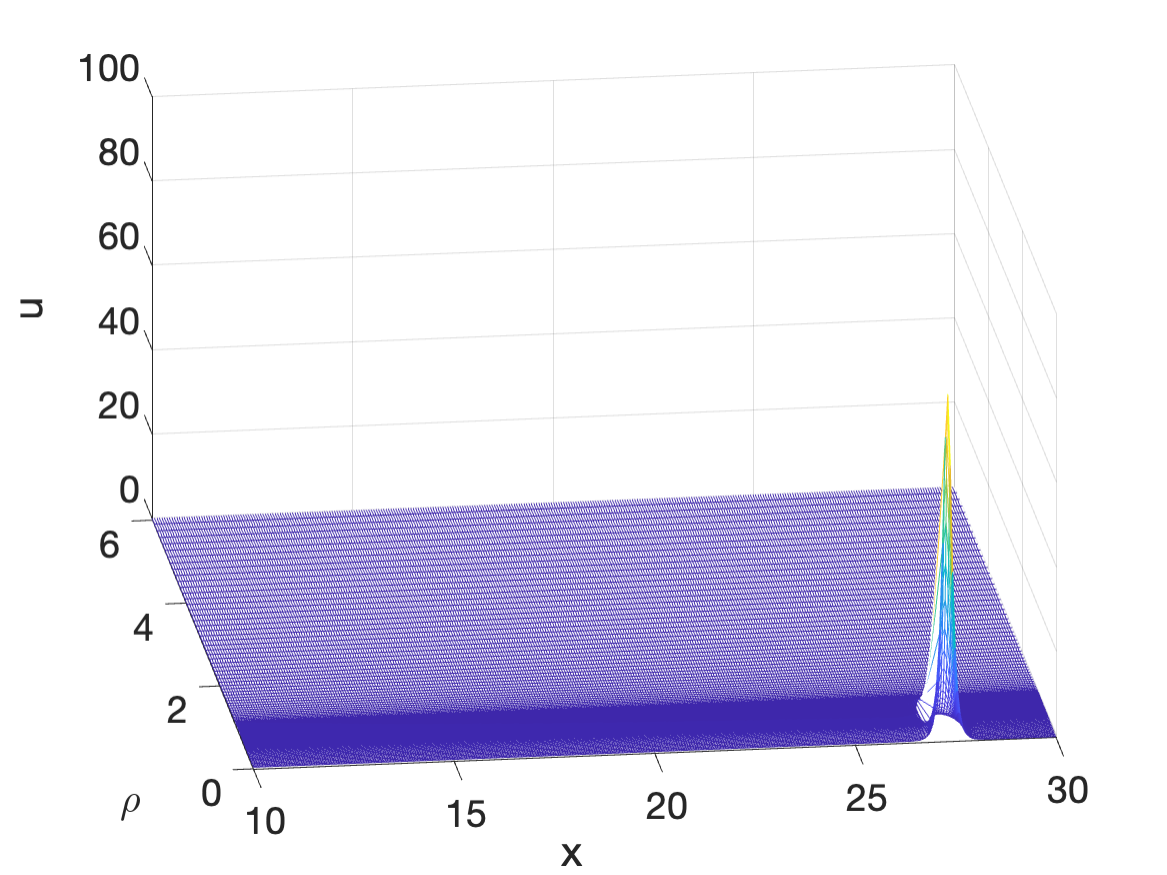}
 \includegraphics[width=0.48\textwidth]{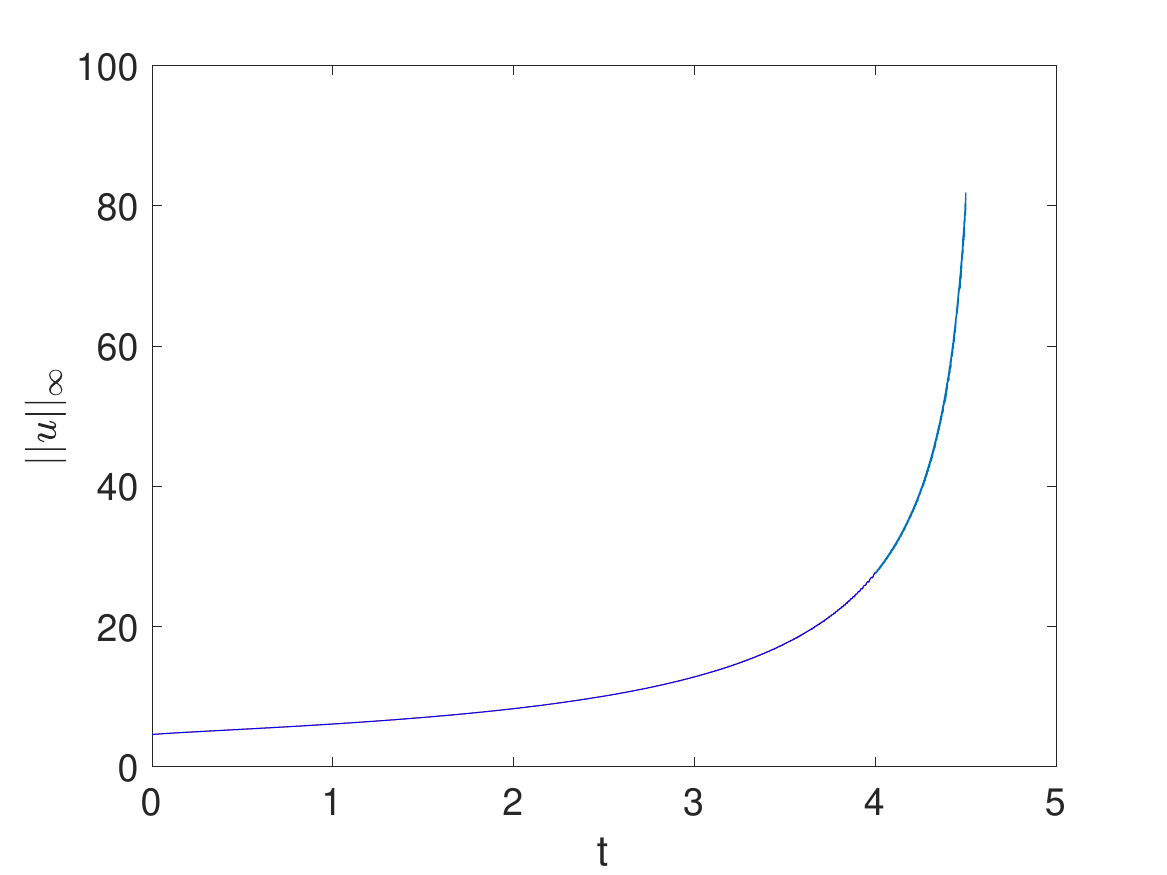}
 \caption{Solution to the 3D critical ZK equation with initial condition 
 $u(0,x,\rho)=1.1\, Q(x,\rho)$, on the left the solution at the final 
 computational time $t= 4.5$, 
 on the right the time dependence of the $L^{\infty}$ norm.} 
 \label{sol1.1}
\end{figure}

For the case of $\lambda=1.1$, we consider the same computational 
domain with $N=2^{12}$ and $N_{t}=10^{3}$ time steps until $t=4$ 
followed by another $10^{3}$ time steps to $t=4.5$. The solution at 
the final time $t=4.5$ can be seen in a close-up on the left of 
Fig.~\ref{sol1.1}. The $L^{\infty}$ norm of the solution on the right 
of the same figure clearly indicates a blow-up in finite time. 

Therefore, it is plausible to conclude that indeed the 3D critical ZK equation when considered on the whole infinite domain 
indeed has the ground state mass as the threshold for the globally existing solutions vs. finite time blow-up solutions, when perturbations of the ground state 
are considered. We next investigate the behavior near this threshold for other type of data, namely, Gaussian, which is a single-peak data with more flattening at the origin and a slightly faster decay at infinity then that of the ground state.


\section{Evolution of Gaussian data}\label{S:gaussian}

In this section we apply our method to study the time evolution of Gaussian initial data, 
\begin{equation}
	u(x,\rho,0) = \lambda e^{-\alpha(x^{2}+\rho^{2})},\quad \lambda >0.
	\label{inigauss}
\end{equation}

Note that the mass of this Gaussian is $M[u] = \lambda^2 \, \big(\frac{\pi}{2\alpha} \big)^{3/2}$. 
Hence, $M[u_0] < \|Q\|_{L^2}^2 \approx 63.78$, when $\lambda < \big(\frac{2\alpha}{\pi} \big)^{3/2} \cdot 63.78$. 
If $\alpha = 1$, then the value on the right-hand side is approximately $5.7$. 
Thus, we expect Gaussian data to disperse if $\lambda \lesssim 5.7$. 
If the sign is reversed, i.e., $M[u_0] > \|Q\|_{L^2}^2$ or $\lambda \gtrsim 5.7$, then we expect the Gaussian initial data to blow up.

Concretely,  we consider the cases:

\qquad (i) $\lambda=5$, which leads to the mass $M[u] \approx 49.2175$, and 

\qquad (ii) $\lambda=6.5$ with the mass $M[u] \sim 83.1776$. 

\begin{figure}[htb!]
 \includegraphics[width=0.5\textwidth]{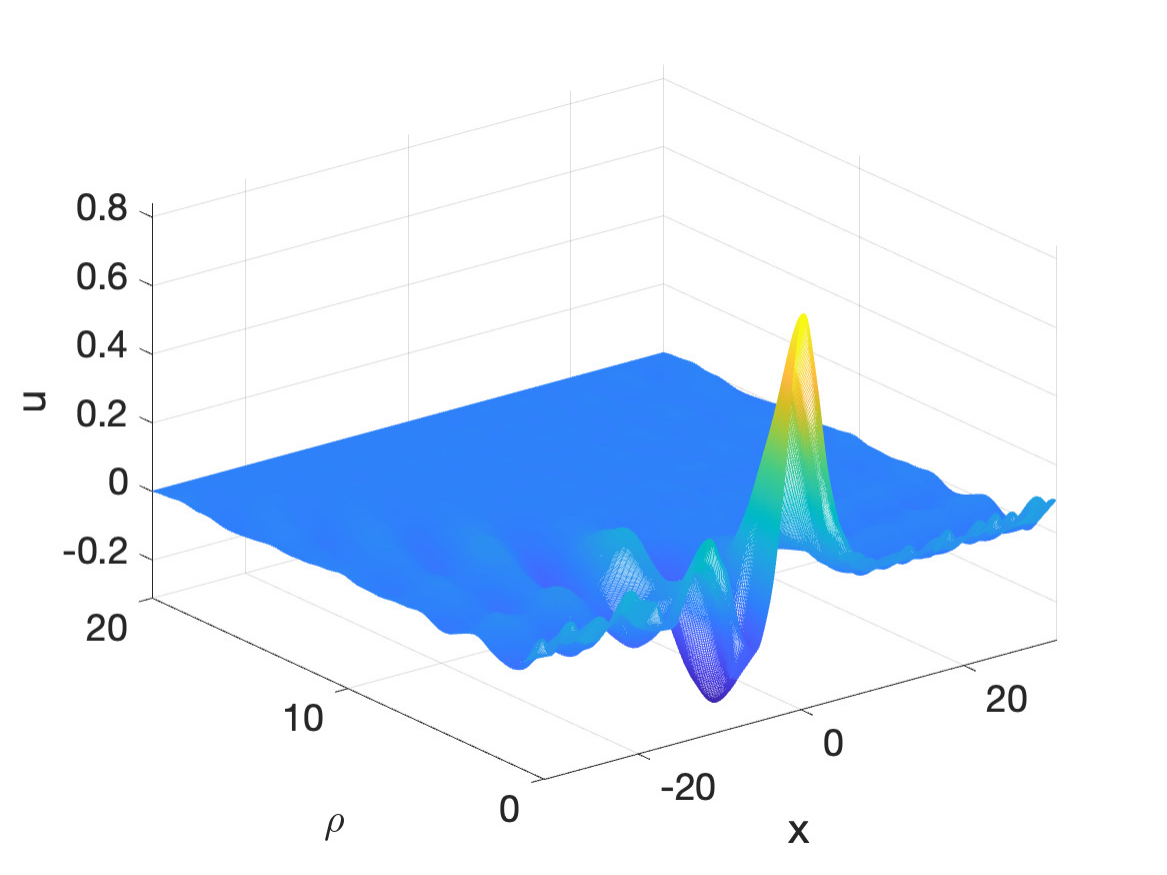}
 \includegraphics[width=0.48\textwidth]{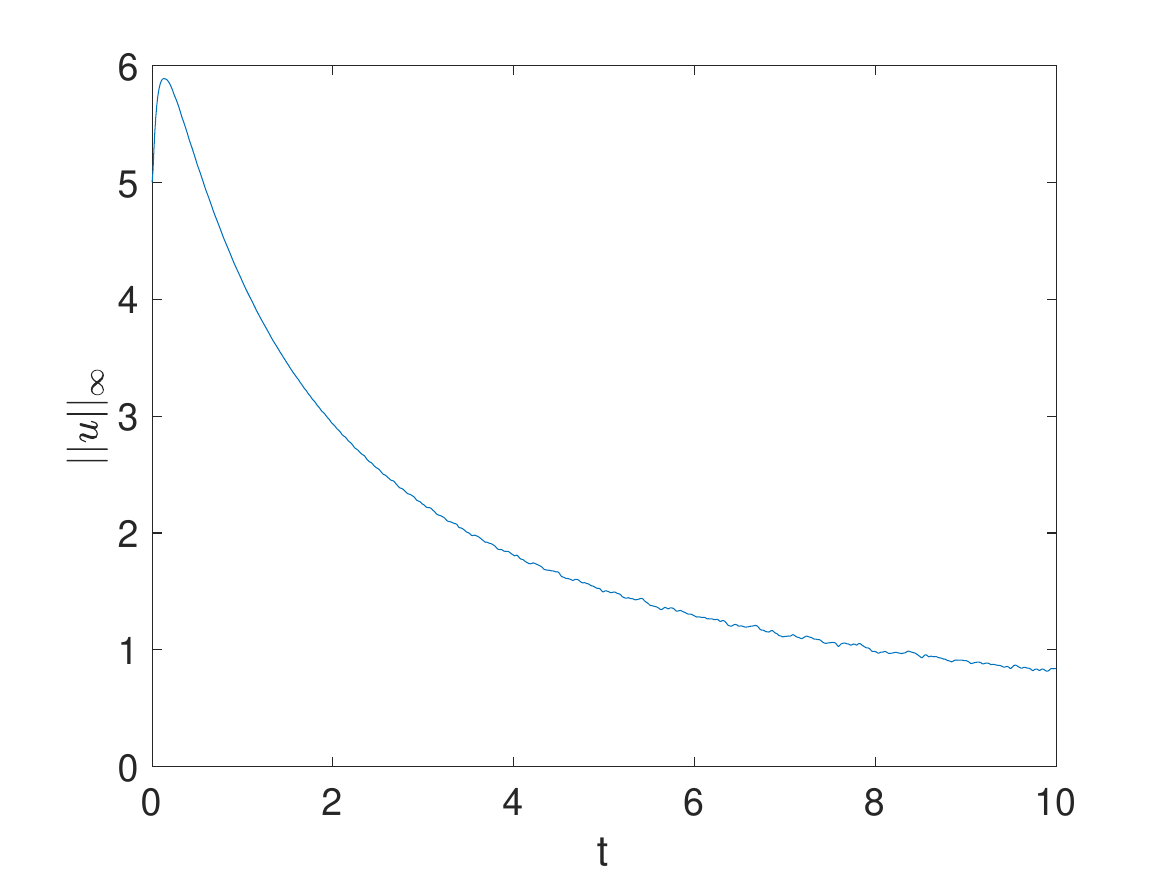}
 \caption{Solution to the 3D critical ZK equation with initial condition  
 $u(0, x,\rho)=5e^{-(x^{2}+\rho^{2})}$: on the left the solution at $t=10$, 
 on the right the time dependence of the $L^{\infty}$ norm.}
 \label{fig5gauss}
\end{figure}

For the case $\lambda=5$ we apply the same numerical parameters as 
in the previous section except for $N_{x}=2^{10}$ and $t\leq 10$. The solution  
clearly disperses as can be seen on Fig.~\ref{fig5gauss} on the left 
where the solution is shown for $t=10$. The initial hump has widened 
and considerably lost height. There is strong radiation in the system 
in the direction of the negative $x$-axis, which re-entered the 
computational domain due to the imposed periodicity in $x$. The $L^{\infty}$ 
norm on the right of the same figure is monotonically decreasing for large $x$ (after a short increase). The 
initial data with smaller mass than the ground state 
appears to be 
simply radiated towards infinity. 

\begin{figure}[htb!]
 \includegraphics[width=0.5\textwidth]{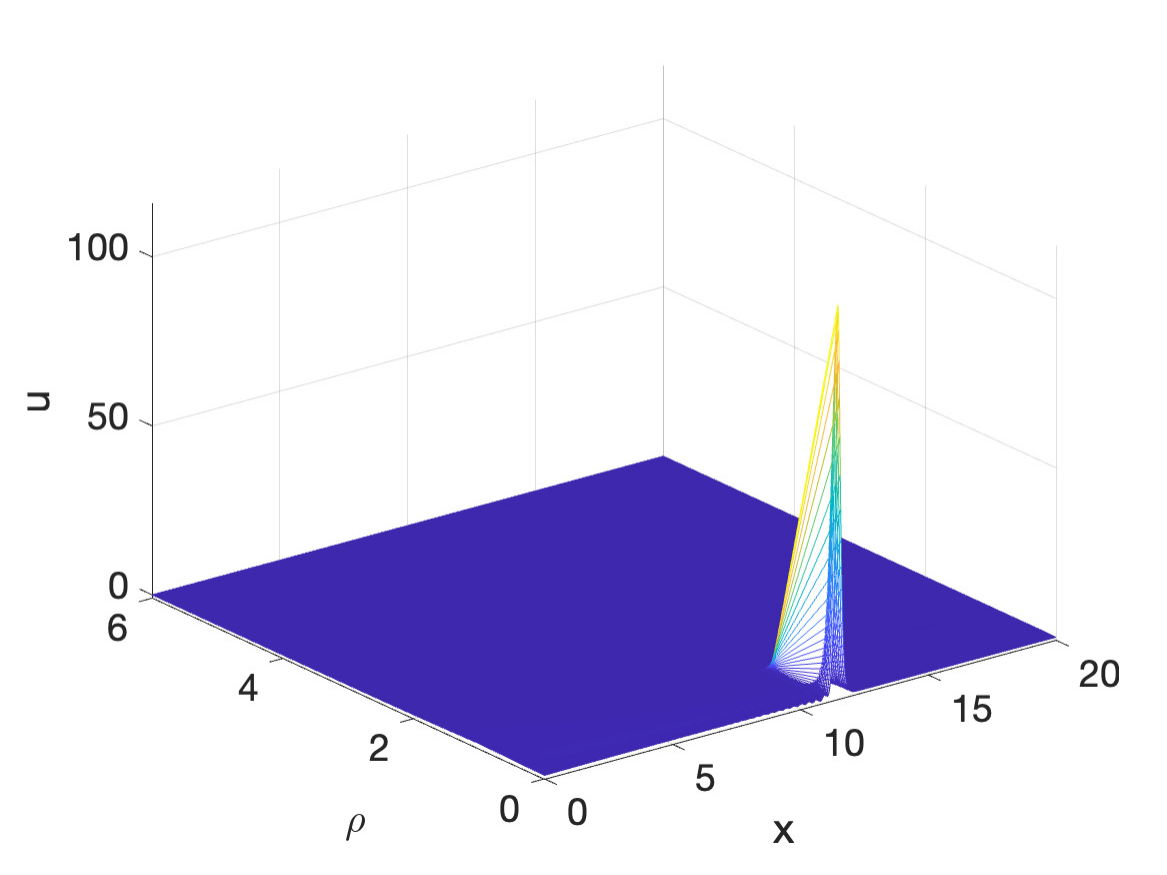}
 \includegraphics[width=0.48\textwidth]{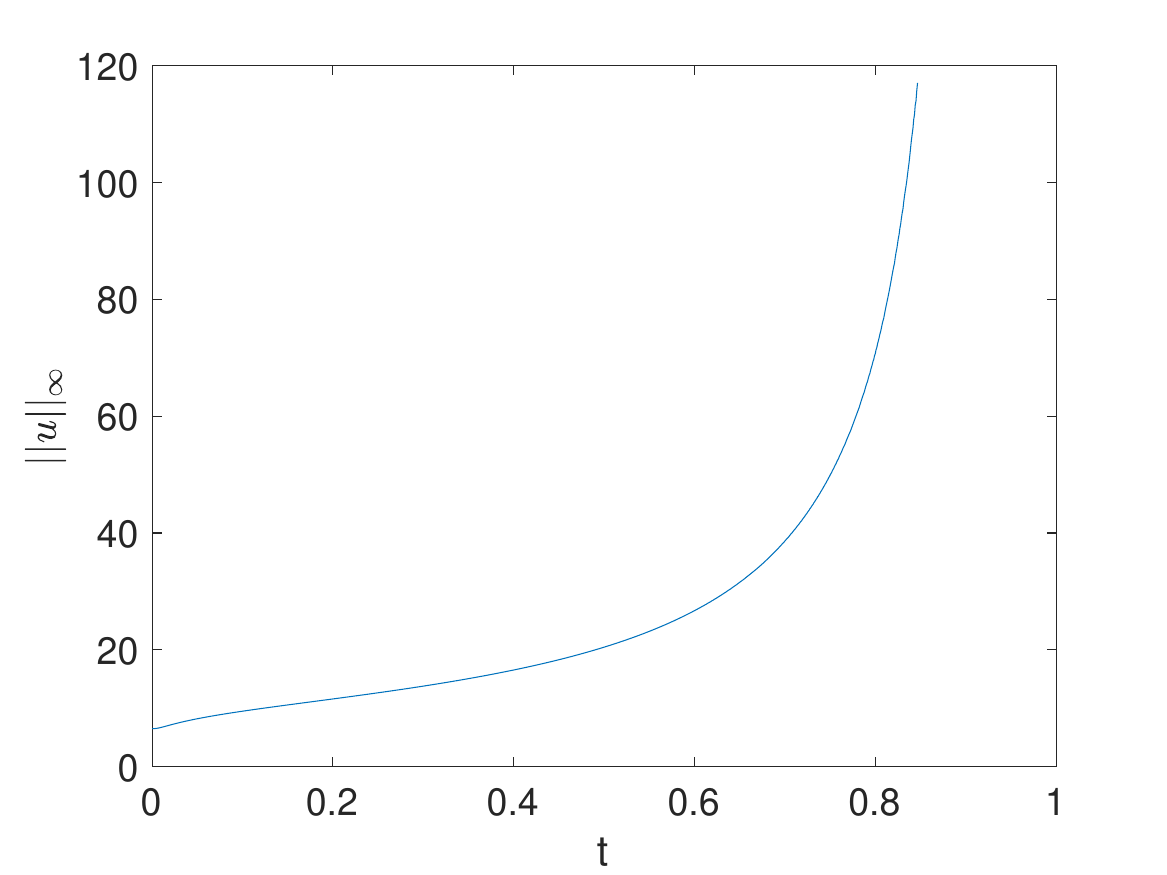}
 \caption{Solution to the 3D critical ZK equation with initial condition  
 $u(0, x,\rho)=6.5e^{-(x^{2}+\rho^{2})}$: on the left the solution at 
 $t=0.85$, 
 on the right the time dependence of the $L^{\infty}$ norm.}
 \label{fig6.5gauss}
\end{figure}
For the case of $\lambda=6.5$, we apply the same numerical parameters 
as for Fig.~\ref{sol1.1}, with $10^{3}$ time steps for the intervals 
$t\in[0,0.75]$ and $t\in[0.75,0.85]$ each. 
The solution at the final 
time can be seen in Fig.~\ref{fig6.5gauss} on the left. The $L^{\infty}$ norm on the right of the same figure indicates once more a blow-up in finite time. 

\begin{figure}[htb!]
\includegraphics[width=0.51\textwidth]{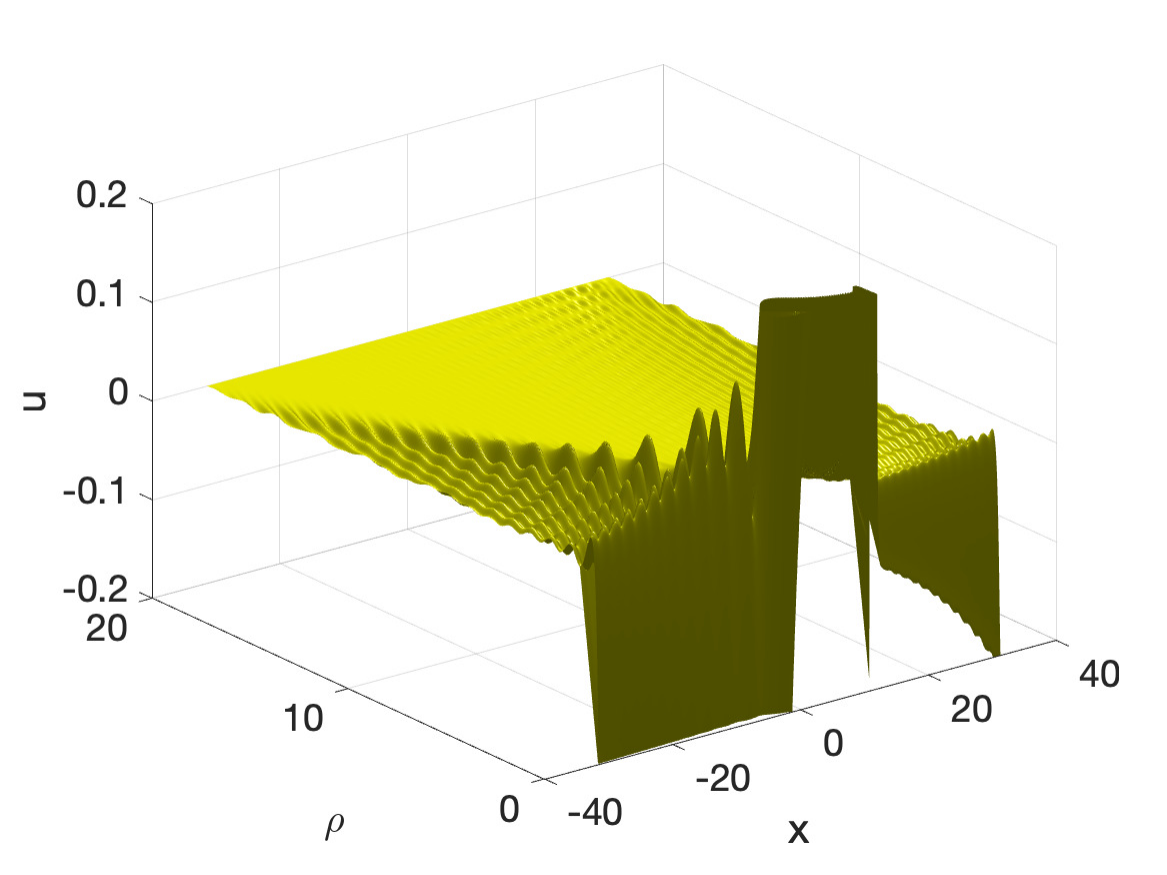}
\includegraphics[width=0.48\textwidth]{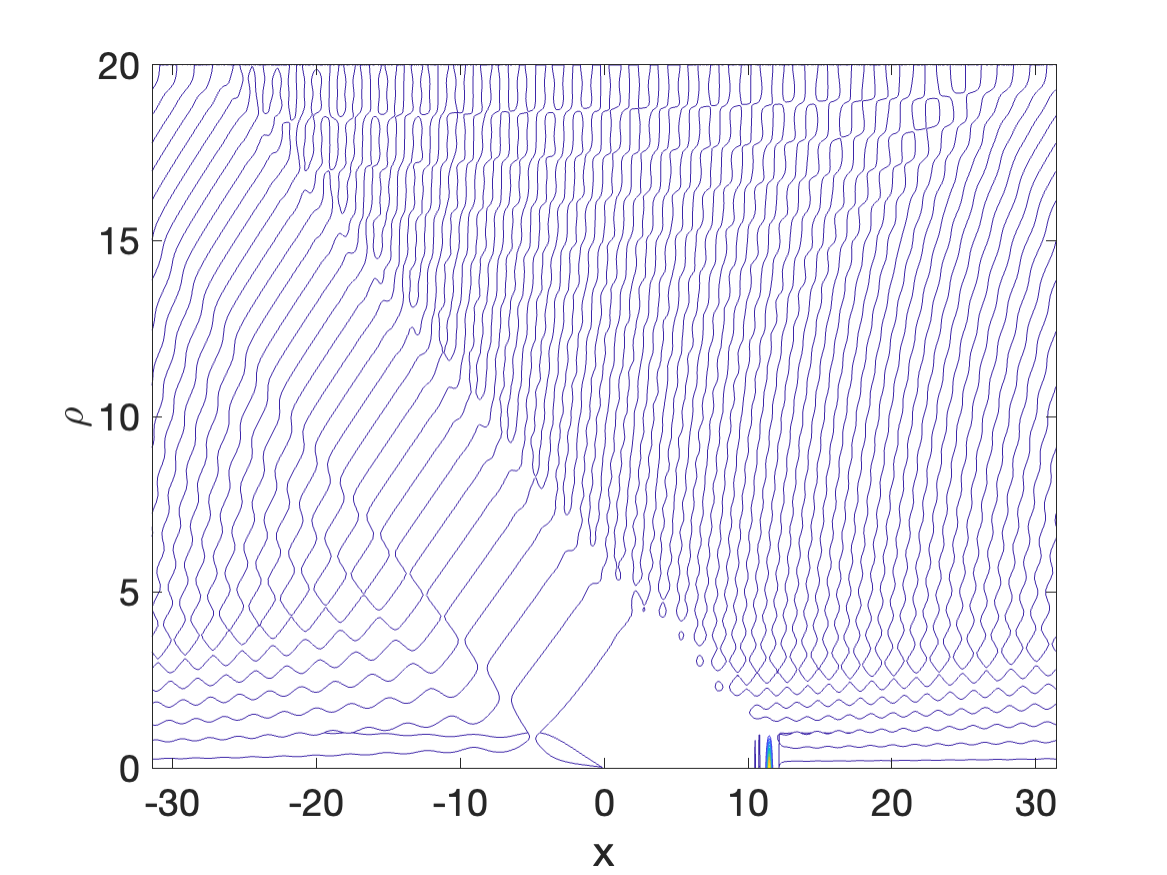}
\caption{Close-up of solution from Fig.~\ref{fig6.5gauss} at the soliton core (left) and contour plot of radiation at final computational time $t=0.85$.}
\label{solfigcontour}
\end{figure}

Finally, in Fig.~ \ref{solfigcontour} the close-up of the solution at the soliton core up to the magnitude of $0.2$ from Fig.~\ref{fig6.5gauss} on the left and the contour plot on the right are shown at the final computational time $t=0.85$ to illustrate spreading of radiation in the previous blow-up example. 
The contour plot of radiation is indicating that the opening of the 
cone around negative $x$-axis has approximately $\tan \theta \approx 
\frac{20}{35}$, which gives half of the radiation opening angle approximately $29.7^o$ degrees, consistent with the $60^o$ radiation angle in ZK equations (e.g., see \cite{FHRY, RRY}). 

\section{Outlook}
In this paper we have presented a numerical approach to the 
higher-dimensional nonlinear dispersive equation with cylindrical 
symmetry on the example of the 3D ZK equation, which has cylindrical symmetry in the $yz$-plane. 
The code uses a Fourier spectral method in $x$ and a multi-domain 
approach in $\rho=\sqrt{y^{2}+z^{2}}$. The domain including the 
origin used $\rho^{2}$ as the independent variable, the second domain used 
$\rho$ with a vanishing condition for the solution at the outer 
boundary (any number of such domains can be included, we concentrate 
here on the simplest case of a single such domain). At the domain 
boundaries the solution is $C^{1}$. The time integration is done with 
an implicit fourth order Runge-Kutta method.

In this paper we concentrated on the $L^{2}$-critical case that is 
most challenging, since a blow-up of initial data is possible in 
finite time as was shown for some examples. The code is, however, 
able to treat general nonlinearities. The type of the blow-up  as 
well as different nonlinearities will be investigated in the future. 

It would be interesting to investigate a non-cylindrically symmetric 
setting, which will need, however, more powerful computational 
architectures to investigate the blow-up solutions. The current code, 
which is essentially 2D because of the cylindrical symmetry,
can be used as a benchmark to develop a full 3D code. 

\bigskip

{\bf Conflict of Interest:} The authors declare that they have no
conflicts of interest.
\bigskip
 
{\bf Data Availability}

Data sharing is not applicable to this article, as no data sets were generated or analyzed.

\end{document}